\documentclass[11pt]{article}
\usepackage{amsmath}
\usepackage{amsmath,amssymb,latexsym,color}
\usepackage[mathscr]{eucal}
\newtheorem{thm}{Theorem}[section]
\newtheorem{defin}[thm]{Definition}
\newtheorem{prop}[thm]{Proposition}

\newtheorem{example}{Example}[section]

\newtheorem{remark}{Remark}[section]

\newcommand{\qed}{\hfill\Box}
\usepackage{CJK}
\begin{document}
\begin{CJK*}{GBK}{song}
\renewcommand{\abovewithdelims}[2]{
\genfrac{[}{]}{0pt}{}{#1}{#2}}

\title{\bf On the metric dimension of bilinear forms graphs}

\author{Min Feng\quad Kaishun Wang\footnote{Corresponding author. E-mail address: wangks@bnu.edu.cn}\\
{\footnotesize   \em  Sch. Math. Sci. {\rm \&} Lab. Math. Com. Sys.,
Beijing Normal University, Beijing, 100875,  China} }
 \date{}
 \maketitle

\begin{abstract}
The metric dimension of a graph is the least number of vertices in a
set with the property that the list of distances from any vertex to
those in the set uniquely identifies that vertex. Bailey and Meagher
\cite{RK} obtained an upper bound on the metric dimension of
Grassmann graphs. In this paper we obtain an upper bound on the
metric dimension of bilinear forms graphs.


\end{abstract}

\bigskip

\bigskip

\section{Introduction}

\bigskip

Let $\Gamma$ be a connected graph. For any two vertices $u$ and $v$,
$d(u,v)$ denotes the distance between $u$ and $v$. By an ordered set
of vertices, we mean a set $W=\{w_{1},\ldots,w_{k}\}$ on which the
ordering $(w_{1},\ldots,w_{k})$ has been imposed. For an ordered
subset $W=\{w_{1},\ldots,w_{k}\}$, we refer to the $k$-vector
$\mathcal{D}(v|W)=(d(v,w_{1}),\ldots,d(v,w_{k}))$ as the
\emph{metric representation of $v$ with respect to $W$}. A
\emph{resolving set} of a graph $\Gamma$ is an ordered subset of
vertices $W$ such that $\mathcal{D}(u|W)=\mathcal{D}(v|W)$ if and
only if $u=v$.
 The \emph{metric dimension} of $\Gamma$, denoted by $\mu(\Gamma)$, is
the smallest size of all the resolving sets of $\Gamma$.

Metric dimension was first introduced in the 1970s, independently by
Harary and Melter \cite{F} and by Slater \cite{PJ}. It is a
parameter that has appeared in various applications, as diverse as
combinatorial optimisation, pharmaceutial chemistry, robot
navigation and sonar. In recent years, a considerable literature has
been developed in graph theory. An interesting case is that of
 distance-regular graphs. For    Johnson graphs and
 Hamming graphs, various results on the metric dimension
have been obtained in \cite{RP,JC,VC,BL,AS}. Recently, Bailey and
Meagher \cite{RK} obtained an upper bound on the metric dimension of
Grassmann graphs. In this paper we consider bilinear forms graphs,
and obtain an upper bound of their matric dimension.

Let $\mathbb{F}_{q}$ be a finite field with $q$ elements. Throughout
this paper, $\mathbb{F}_{q}^{n+d}$ denotes the $(n+d)$-dimensional
vector space over $\mathbb{F}_{q}$, and $N$ denotes a fixed
$n$-dimensional subspace of $\mathbb{F}_{q}^{n+d}$.

The \emph{bilinear forms graph} $H_{q}(n,d)$ has  as its vertex set
the set of all $d$-dimensional subspaces of $\mathbb{F}_{q}^{n+d}$
intersecting trivially with $N$, and two vertices are adjacent if
they intersect in a subspace of dimension $d-1.$ The bilinear forms
graph $H_{q}(n,d)$ is a distance-regular graph with $q^{nd}$
vertices and diameter $\min(n,d)$ such that the distance between two
vertices $A$ and $B$ is $d-\dim(A\cap B)$. For more information
about distance-regular graphs, we refer readers to \cite{AE}.

 Note that $H_{q}(n,1)$ is a complete graph
  whose metric dimension   is $q^{n}-1$. Also,
$H_{q}(n,d)$ is isomorphic to $H_{q}(d,n)$, we only need to consider
the case $n\geq d\geq 2$.

In this paper, we obtain the following   result:

\begin{thm}\label{dl}
Let $n\geq d\geq2$. Then   the metric dimension   of the bilinear
forms graph $H_{q}(n,d)$ satisfies
\[
\mu(H_{q}(n,d))\leq\left\{
\begin{array}{cc}
q^{n+d-1}  &\textup{if}~n\geq d+2,\\
q^{n+d}~~~ &\textup{otherwise}.
\end{array}\right.
\]
\end{thm}

\section{Proof of Theorem \ref{dl}}

We shall  prove Theorem \ref{dl} by constructing resolving sets. Our
construction   requires some notion from finite geometry.

A {\em partition} of the vector space \emph{V} is a set
$\mathcal{P}$ of subspaces of $V$ such that any non-zero vector is
contained in exactly one element of $\mathcal{P}$. If $T=\{ \dim
W\mid W\in \mathcal{P}\}$, the partition $\mathcal{P}$ is said to be
a $T$-$partition$ of $V$.

\begin{prop}\label{pro}
(\cite[Lemma 2]{AB}) Let $s$ and $t$ be positive integers with
$s+t=n+d$, then there exists an $\{s,t\}$-partition of
$\mathbb{F}_{q}^{n+d}$.
\end{prop}

\noindent {\em Proof of Theorem \ref{dl}}: We divide the proof in
two cases:

\medskip
{\em Case 1}. $n\geq d+2$.
 The proof of Proposition \ref{pro} implies that
$\mathbb{F}_{q}^{n+d}$ has an $\{n-1,d+1\}$-partition
$\mathcal{P}_{1}=\{\tilde{N},W_{1},\ldots,W_{m}\}$, where
$\tilde{N}\subset N,~\dim \tilde{N}=n-1$ and $\dim
W_{i}=d+1,~i=1,\ldots,m$. Note that  $m=q^{n-1}$.

For each $i$, let $N_{i}=W_{i}\cap N$. Since $\dim(W_{i}+N)=n+d$,
$\dim N_{i}=1$.   Suppose
$$
\mathcal{M}=\bigcup_{i=1}^{m} \mathcal{M}_i,
$$
where $\mathcal{M}_i$ is the collection of   $d$-subspace of $W_{i}$
intersecting trivially with $N_{i}$.
  For any $U\in \mathcal{M}$,
$U+N=\mathbb{F}_{q}^{n+d}$, so $U\cap N=\{0\}$. It follows that
$\mathcal{M}$ is a subset of the vertex set of $H_{q}(n,d)$.

Next we shall prove $\mathcal{M}$ is a resolving set of
$H_{q}(n,d)$. We only need to show that, for any two distinct
vertices, there exists a vertex $U\in\mathcal{M}$ such that
\begin{equation}\label{inq}
\dim(A\cap U)\neq \dim(B\cap U).
\end{equation}
For each $i$, let $A_{i}=A\cap W_{i}$, $B_{i}=B\cap W_{i}$. Since
$A\neq B$, there exists an $i$ such that $A_{i}\neq B_{i}$. Suppose
$\dim A_{i}=s\leq\dim B_{i}=t$. Let $\{\beta_{1},\ldots,\beta_{t}\}$
be a basis for $B_{i}$ and $\{\theta\}$ be a basis for $N_{i}$.

\medskip
 {\em Case 1.1.}  $s<t$. Note that $\{\theta,\beta_{1},
\ldots,\beta_{t}\}$ is linearly independent. Extend this to a basis
$\{\theta,\beta_{1},\ldots,\beta_{t},\gamma_{1},\ldots,\gamma_{d-t}\}$
for $W_{i}$; let $U$ be the $d$-dimensional space spanned by
$\{\beta_{1},\ldots,\beta_{t},\gamma_{1},\ldots,\gamma_{d-t}\}$. By
construction, $U$ is an element of $\mathcal{M}$ satisfying
$$
\dim(A\cap U)=\dim(A_{i}\cap U)\leq \dim A_{i}=s<t=\dim
B_{i}=\dim(B_{i}\cap U)= \dim(B\cap U),
$$
 so (\ref{inq}) holds.

\medskip
 {\em Case 1.2.} $s=t$.  Since $A_{i}\neq B_{i}$, there exists an $\alpha\in
A_{i}\backslash B_{i}$. Then each of $\{\alpha,\beta_{1},
\ldots,\beta_{t}\}$, $\{\alpha,\theta\}$  and $\{\theta,\beta_{1},
\ldots,\beta_{t}\}$ is linearly independent.

\medskip
 {\em Case 1.2.1.}   $\{\alpha,\theta,\beta_{1},
\ldots,\beta_{t}\}$ is linearly dependent. Extend
$\{\theta,\beta_{1}, \ldots,\beta_{t}\}$ to a basis
$\{\theta,\beta_{1},\ldots,\beta_{t},\gamma_{1},\ldots,\gamma_{d-t}\}$
for $W_{i}$ and let $U$ be the $d$-dimensional space spanned by
$\{\beta_{1},\ldots,\beta_{t},\gamma_{1},\ldots,\gamma_{d-t}\}$.
Since   $\alpha\not\in U$, $U$ is an element of $\mathcal{M}$ such
that
$$
\dim(A\cap U)=\dim(A_{i}\cap U)< \dim A_{i}=\dim(B\cap U),
$$
so (\ref{inq}) holds.

\medskip
 {\em Case 1.2.2.}    $\{\alpha,\theta,\beta_{1},
\ldots,\beta_{t}\}$ is linearly independent. Extend this to the
basis
$\{\alpha,\theta,\beta_{1},\ldots,\beta_{t},\gamma_{1},\ldots,\gamma_{d-t-1}\}$
for $W_{i}$ and let $U$ be the $d$-dimensional space spanned by
$\{\alpha+\theta,\beta_{1},\ldots,\beta_{t},\gamma_{1},\ldots,\gamma_{d-t-1}\}$.
Since   both $\{\theta,\alpha+\theta,\beta_{1},
\ldots,\beta_{t},\gamma_{1},\ldots,\gamma_{d-t-1}\}$ and
$\{\alpha,\alpha+\theta,\beta_{1},
\ldots,\beta_{t},\gamma_{1},\ldots,\gamma_{d-t-1}\}$ are linearly
independent, we have $\theta\not\in U$ and $\alpha\not\in U$.
Consequently, there exists a $U\in\mathcal{M}$ such that
$$
\dim(A\cap U)< \dim(B\cap U),
$$
so   (\ref{inq}) holds.

\medskip

By \cite[Lemma 9.3.2]{AE}, $|\mathcal{M}|=q^{n+d-1}$. Hence,
$\mu(H_{q}(n,d))\leq q^{n+d-1}.$

\medskip{\em Case 2}. $d\leq n\leq d+1$.

By the proof of Proposition \ref{pro}, $\mathbb{F}_{q}^{n+d}$ has an
$\{n,d\}$-partition $\mathcal{P}_{2}=\{N,V_{1},\ldots,V_{m}\}$,
where $\dim V_{i}=d,~i=1,\ldots,m$. Note that $m=q^{n}$.

Let $\bar{N}$ be a fixed $1$-dimensional subspace of $N$. For each
$i$, let $W_{i}=\bar{N}+V_{i}$. Suppose
$$
\mathcal{M}=\bigcup_{i=1}^{m} \mathcal{M}_i,
$$
where $\mathcal{M}_i$ is the collection of   $d$-subspace of $W_{i}$
intersecting trivially with $\bar N$.  For any $U\in \mathcal{M}$,
$U+N=\mathbb{F}_{q}^{n+d}$, so $U\cap N=\{0\}$. It follows that
$\mathcal{M}$ is a subset of the vertex set of $H_{q}(n,d)$.

Similar to Case 1,   $\mathcal{M}$ is a resolving set of
$H_{q}(n,d)$ with $|\mathcal{M}|=q^{n+d}$. Hence,
$\mu(H_{q}(n,d))\leq q^{n+d}.$

\medskip

By above discussion, we complete the proof. $\qed$

Babai \cite{LB} obtained bounds on a parameter of primitive
distance-regular graphs which is equivalent to the metric dimension.
A natural question is to compare our result with those. For the case
of the bilinear forms graph $H_{q}(n,d)$, Babai's most general bound
(see \cite[Theorem 2.1]{LB}) yields
\begin{eqnarray*}\label{b1}
\mu(H_{q}(n,d))<4\sqrt{q^{nd}}\log(q^{nd}),
\end{eqnarray*}
while his stronger bound (see \cite[Theorem 2.4]{LB}) yields
\begin{eqnarray*}\label{b2}
\mu(H_{q}(n,d))<2d\cdot\frac{q^{nd}}{q^{nd}-M}\cdot\log(q^{nd}),
\end{eqnarray*}
where
$$
M=\max_{0\leq i\leq
d}\abovewithdelims{n}{i}_{q}\abovewithdelims{d}{i}_{q}(q^{i}-1)(q^{i}-q)\cdots(q^{i}-q^{i-1}).
$$

For $n, d \geq 4$,  our   bound is better than Babai's most general
bound. For the left case, there is some $q$ such that our bound is
better. Babai's stronger bound   is difficult to evaluate exactly,
so we conduct   some experiments using MATLAB to compare this bound
with our bound. We find our bound is better in most cases for $q=2$.

\section*{Acknowledgement} The authors would like to thank  Yuefeng Yang for communicating.   This
research is supported by NSF of China (10871027) and NCET-08-0052.

\end{CJK*}

\end{document}